\documentclass{amlts}
\begin{document}
\annalsline{155}{2002}
\received{May 12, 2000; corrected version November 11, 2026}
\startingpage{789}
 \font\tenrm=cmr10
\input amssym.def
\input amssym.tex
\def\eqref#1{(\ref{#1})}
\catcode`\@=11
\font\twelvemsb=msbm10 scaled 1100
\font\tenmsb=msbm10
\font\ninemsb=msbm10 scaled 800
\newfam\msbfam
\textfont\msbfam=\twelvemsb  \scriptfont\msbfam=\ninemsb
  \scriptscriptfont\msbfam=\ninemsb
\def\msb@{\hexnumber@\msbfam}
\def\Bbb{\relax\ifmmode\let\next\Bbb@\else
 \def\next{\errmessage{Use \string\Bbb\space only in math
mode}}\fi\next}
\def\Bbb@#1{{\Bbb@@{#1}}}
\def\Bbb@@#1{\fam\msbfam#1}
\catcode`\@=12

 \catcode`\@=11
\font\twelveeuf=eufm10 scaled 1100
\font\teneuf=eufm10
\font\nineeuf=eufm7 scaled 1100
\newfam\euffam
\textfont\euffam=\twelveeuf  \scriptfont\euffam=\teneuf
  \scriptscriptfont\euffam=\nineeuf
\def\euf@{\hexnumber@\euffam}
\def\frak{\relax\ifmmode\let\next\frak@\else
 \def\next{\errmessage{Use \string\frak\space only in math
mode}}\fi\next}
\def\frak@#1{{\frak@@{#1}}}
\def\frak@@#1{\fam\euffam#1}
\catcode`\@=12

\newcommand{\C}{{\Bbb C}}
\newcommand{\R}{{\Bbb R}}
\newcommand{\B}{{\mathcal B}}
\newcommand{\Z}{{\Bbb Z}}
\newcommand{\BMO}{\mathop{\rm {\rm BMO}}\nolimits}
\renewcommand{\Re}{\mathop{\rm Re}\nolimits}
\renewcommand{\Im}{\mathop{\rm Im}\nolimits}
\newcommand{\dist}{\mathop{\rm dist}\nolimits}
\newcommand{\eps}{\varepsilon}
\newcommand{\dbar}{\bar\partial}

\title{Fourier frames}
  \acknowledgements{The first author is supported by the DGICYT grant
PB98-1242-C02-01 and by the DURSI grant 2001SGR-00172. Part of this work
was done June--July 1999. The second author participated in and 
was supported by the
Semester on Analysis  at Centre de Recerca Matem\`atica,
Barcelona, 1999. This is a corrected version of the paper published in Ann.\ of Math.\
(2) \textbf{155} (2002), 789--806. The typos and mistakes corrected here were found during a
formalization of the paper in Lean~4 (\texttt{github.com/joaquimortega/MR1923965}); the changes
are kept to a minimum.}
  \twoauthors{Joaquim Ortega-Cerd\`a}{Kristian Seip}
 \institutions{Universitat de Barcelona,  Barcelona, Spain
\\
 {\eightpoint {\it E-mail address\/}: quim@mat.ub.es}\\
\vglue6pt
Norwegian University of Science and Technology, Trondheim, Norway
\\
 {\eightpoint {\it E-mail address\/}: seip@math.ntnu.no
}}
\vglue24pt
\centerline{\bf Abstract}
\vglue12pt

We solve the problem of Duffin and Schaeffer (1952) of
characterizing those sequences of real frequencies which generate
Fourier frames. Equivalently, we characterize the sampling
sequences for the Paley-Wiener space. The key step is to connect
the problem with de Branges' theory of Hilbert spaces of entire
functions. We show that our description of sampling sequences
permits us to obtain a classical inequality of H.~Landau as a
consequence of Pavlov's description of Riesz bases of complex
exponentials and the John-Nirenberg theorem. Finally, we discuss
how to transform our description into a working condition by
relating it to an approximation problem for subharmonic functions.
By this approach, we determine the critical growth rate of a
nondecreasing function $\psi$ such that the sequence
$\{\lambda_k\}_{k\in\Z}$ defined by $\lambda_k+\psi(\lambda_k)=k$
is sampling.

\section{Introduction}

Following Duffin and Schaeffer \cite{DS52}, we say that a
system of complex exponentials $\{e^{i\lambda_k x}\}$, with
$\Lambda=\{\lambda_k\}_{k\in \Z}$ a sequence of real numbers, is a
{\it Fourier frame} if there exist positive constants $A$ and $B$
such that
$$A \int_{-\pi}^\pi |f(x)|^2 dx \leq
\sum_{k=-\infty}^\infty \left|\int_{-\pi}^\pi f(x)e^{-i\lambda_k x}
dx\right|^2 \leq B \int_{-\pi}^\pi |f(x)|^2 dx $$ for all $f\in
L^2(-\pi,\pi)$. The purpose of this paper is to give a
description of those sequences $\Lambda$ that generate Fourier frames.\pagebreak

The main result of \cite{DS52} is a sufficient density
condition for $\{e^{i\lambda_k x}\}$ to constitute a Fourier
frame. Over the past two decades, the work of Duffin and Schaeffer
has been remarkably influential, but mainly so because of its
description of abstract frame expansions as an alternative to
orthonormal bases; see, e.g., \cite{Da92}. There has been little
progress on the problem of characterizing the original Fourier
frames of Duffin and Schaeffer.

Our work relies on several profound results. Specifically, it may be seen as
an interplay between three themes:

\demo{{\rm 1)} de Branges\/{\rm '} Hilbert spaces of entire functions} Our main
theorem is based on de Branges' theory \cite{dB68}. The Hilbert space
of interest to us is the classical Paley-Wiener space,
the prime example of a de Branges space. We denote this space by ${\rm PW}$;
it consists of all entire functions of exponential type at
most $\pi$ whose restrictions to $\R$ are square-integrable. We say that
a sequence of real numbers $\Lambda=\{\lambda_k\}_{k\in\Z}$ is
{\it sampling} for ${\rm PW}$ if there exist positive constants $A$ and $B$ such
that
$$
A\int_{-\infty}^{\infty}|f(x)|^2\, dx\leq
\sum_{k=-\infty}^{\infty}|f(\lambda_k)|^2
\leq B \int_{-\infty}^{\infty}|f(x)|^2\, dx
$$
holds for all $f\in {\rm PW}$. By the Paley-Wiener theorem,
{\it $\Lambda$ is sampling for ${\rm PW}$ if and only if
$\{e^{i\lambda_k x}\}$ is a Fourier frame}. In other words, our task is to
describe the sampling sequences for ${\rm PW}$.

The Paley-Wiener space is a Hilbert space when equipped with its
standard norm, the $L^2(\R)$-norm. The symbol ${\rm PW}$ will stand for
this particular Hilbert space. However, if $\Lambda$ is a sampling
sequence, we may equip the Paley-Wiener space with an equivalent
norm given by $\sqrt{\sum_k|f(\lambda_k)|^2}$. In this way, we
obtain a different Hilbert space, which is again a de Branges
space. This simple observation has the amazing consequence that
several basic results of de Branges' theory apply to our problem.
In particular, two results of de Branges about {\it norm
identities} may be combined to yield a result about {\it norm
equivalence} (Theorem 1 below).
\enddemo

The required background on de Branges spaces is reviewed in
Section 2, and then our main theorem --- a necessary and sufficient
condition for sampling --- is stated and proved in Section 3.

\demo{{\rm 2)} Riesz bases of complex exponentials} Fourier frames are
intimately connected with Riesz bases of complex exponentials, or
equivalently, sampling sequences are intimately connected with
complete interpolating sequences (see below for definition).
This is of interest to us, because there exists a beautiful description
of the latter kind of sequences, due to Pavlov \cite{Pa79}; see
\cite{HNP81} for an extensive account and \cite{LS97}
for some recent progress on this topic.
\enddemo

We say that $\Lambda$ is a
{\it complete interpolating sequence} for ${\rm PW}$ if the
interpolation problem $f(\lambda_k)=a_k$ has a unique solution
$f\in {\rm PW}$ whenever $\{a_k\}$ is square-summable, or equivalently,
if it is sampling but fails to be so on the removal of any one of
the points $\lambda_k$. Thus a complete interpolating sequence
is a sampling sequence with no ``redundant'' points. The sequence of
integers is of course the leading example of such a sequence.
\vglue4pt

Some additional terminology is required to state Pavlov's theorem.
Let us assume that $\lambda_k\le \lambda_{k+1}$ for all $k$.
A sequence $\Lambda$ is {\it separated} if
$$q=\inf_k(\lambda_{k+1}-\lambda_k)>0;$$ $q$
is referred to as the {\it separation constant} of $\Lambda$.
With a separated sequence $\Lambda$ we associate a distribution function
$n_\Lambda(t)$ defined such that for $a<b$
$$ n_\Lambda(b)-n_\Lambda(a)=\# (\Lambda\cap (a,b]),
$$
and normalized such that $n_\Lambda(0)=0$. There is clearly a one-to-one
correspondence between $\Lambda$ and $n_\Lambda$. It is
plain that all complete interpolating sequences are  separated.

Let
$$
P_y[h](x)=\frac y\pi\int_{-\infty}^\infty \frac{h(t)}{(x-t)^2+y^2}\, dt
$$
denote the Poisson extension of $h$ to the upper half-plane $y>0$,
and let $u\mapsto \tilde{u}$ stand for the usual conjugation
operator. The following theorem, with the Helson-Szeg\"{o}
condition appearing in slight disguise, is a reformulation (see
\cite[p. 286]{HNP81}) of Pavlov's original result.

\proclaimtitle{Pavlov; Hru\v{s}\v{c}ev, Nikol'ski{\u\i}, Pavlov}
\specialnumber{A}\proclaim{Theorem}
A separated sequence $\Lambda$ of real numbers is a complete interpolating
sequence for ${\rm PW}$ if and only if $h(t)=n_\Lambda(t)-t$
is a function in $\BMO(\R)$ such that
$P_1[h](x)=\tilde{u}(x)+v(x)+C${\rm ,} with
$u,v\in L^\infty(\R)$ and $\|v\|_\infty< 1/4${\rm .}
\endproclaim

To see  the problem of relating this result to sampling sequences,
we begin by stating two density conditions for sampling. The first is a
fairly elementary result: {\it $\Lambda$ is sampling if
there exist positive constants $\varepsilon$ and $C$ such that
$$  n_\Lambda(b)-n_\Lambda(a)\geq (1+\varepsilon)(b-a) -C
$$
for all} $a<b$. This is essentially what was obtained by Duffin
and Schaeffer (see also Theorem 2.1 of \cite{Se95}). A more
profound result is the following inequality\footnote{In
\cite{La67}, the inequality is given in the following more
general form. Let $\Omega$ be a finite union of intervals with
$|\Omega|=2\pi$ and ${\rm PW}(\Omega)$ the subspace of $L^2(\R)$
consisting of functions whose Fourier transforms vanish outside
$\Omega$; define sampling sequences for ${\rm PW}(\Omega)$ as above.
Then the inequality remains true if we replace ${\rm PW}$ by
${\rm PW}(\Omega)$.} of Landau \cite{La67}.

\proclaimtitle{Landau}
\specialnumber{B}\proclaim{Theorem} If $\Lambda$ is a separated sampling sequence
for ${\rm PW}${\rm ,} then
$$
n_\Lambda(b)-n_\Lambda(a)\ge b-a -A \log^+(b-a) -B
$$
when $a<b${\rm ,} with constants $A$ and $B$ independent of $a,b${\rm .}
\endproclaim

Landau's inequality is best possible, as shown by an example in
Section 4. We shall see that Theorem A and our main theorem
(Theorem 1) imply that {\it Landau\/{\rm '}\/s inequality is a consequence
of the John\/{\rm -}\/Nirenberg theorem for $\BMO$ functions.} This result
will be proved in Section 4 below.
\vglue2pt

The two inequalities just recorded have the following
consequence. Set
$$
D^-(\Lambda)=\lim_{R\to\infty}\frac{\min_x (n_\Lambda(x+R)-n_\Lambda(x))}{R},
$$
which is the Beurling lower uniform density of $\Lambda$. Then
$\Lambda$ is sampling if $D^-(\Lambda)>1$ and fails to be sampling
if $D^-(\Lambda)<1$. (See also \cite{Ja91} for a somewhat
different formulation.) It was proved in \cite{Se95} that when
$D^-(\Lambda)>1$, there exists a subsequence
$\Lambda'\subset\Lambda$ such that $\Lambda'$ is a complete
interpolating sequence, but an example was given of a sampling
sequence $\Lambda$ for which no subsequence is a complete
interpolating sequence. {\it The latter type of sampling
sequences are the only ones which are essentially different from
complete interpolating sequences}. The third theme of this paper
is what is needed to explore their subtle properties.

\demo{{\rm 3)} Approximation of subharmonic functions} We have in
mind the following general type of problem. Given a subharmonic
function $U$, find an entire function $f$ such that $\log|f|$ is
in some sense a good approximation to $U$. Our main theorem
becomes a working condition only if we are able to solve certain
problems of this kind. An exhaustive discussion of the solutions
relevant for our sampling condition seems impracticable. Instead,
we use an approach of Lyubarskii and Malinnikova [LM01]
to elaborate an illustrative example. In Section 5, we determine
the critical growth rate of a nondecreasing function $\psi$ such
that $\Lambda=\{\lambda_k\}_{k\in\Z}$ defined by
$\lambda_k+\psi(\lambda_k)=k$ is sampling. Roughly speaking, we
obtain that if $\psi(x)=0$ for $x<0$, $\psi$ is sufficiently
regular, and $\psi(x)\to\infty$ when $x\to\infty$, then $\Lambda$
is sampling if and only if $\psi(x)$ tends to $\infty$ at least as
fast as $\log^+x$.

This result may seem surprising: A sequence is sampling if and only if it
constitutes a sufficiently large deviation from $\Z$. The central point
here is that the sampling condition implies that we need to solve the above
approximation problem for the logarithmic potential of the measure $d\psi$.
Then difficulties arise if the measure is too much ``spread out''.\pagebreak

\section{Preliminaries on de Branges spaces}

A Hilbert space $H$ of entire functions is a {\it de Branges space} if the
following conditions are met:
\begin{itemize}
\item[(H1)] Whenever $f$ is in the space and has a nonreal zero $\zeta$,
the function $g(z)=f(z)(z-\bar{\zeta})/(z-\zeta) $ is in the space
and has the same norm as $f$.
\item[(H2)] For every nonreal $\zeta$ the linear functional
defined on the space by
$f\mapsto f(\zeta)$ is continuous.
\item[(H3)] The function $f^*(z)=\overline{f(\bar{z})}$ belongs to the space
whenever $f$ belongs to the space and $f^*$ has the same norm as $f$.
\end{itemize}
The simplest model example of a de Branges space is ${\rm PW}$. A wider class of
examples
can be constructed via functions from the Hermite-Biehler class of entire
functions. We denote this class by $\overline{\rm HB}$; it consists
of all entire functions $E$ with no zeros in the upper half-plane and
satisfying $|E(z)|\geq |E(\overline{z})|$ whenever $\Im z>0$.
To every function $E\in \overline{\rm HB}$ we associate a Hilbert space
$H(E)$, which consists of all entire functions
$f$ such that both $f(z)/E(z)$ and $f^*(z)/E(z)$ belong
to $H^2$ of the upper half-plane; the $H(E)$-norm of $f$ is given by
$$
\|f\|_E^2=
\int_{-\infty}^\infty \frac{|f(t)|^2}{|E(t)|^2}\, dt.
$$
It is clear that $H(E)$ is a de Branges space. The following fundamental
theorem  of de Branges says that all de Branges spaces arise in this way
\cite[p. 57]{dB68}.

\specialnumber{C}\proclaim{Theorem}
A Hilbert space $H$ whose elements are entire functions{\rm ,} which satisfies
  {\rm (H1),}  {\rm (H2),} and  {\rm (H3),} and which contains a nonzero element{\rm ,}
is equal isometrically to some space $H(E)${\rm .}
\endproclaim

It follows from condition  (H2)  that for each nonreal $\zeta$
there exists a reproducing
kernel $K(\zeta, z)$ for the space   $H(E)$.
The kernel has the following representation, which also
defines it for real $\zeta$ \cite[p. 50]{dB68}.

\specialnumber{D}\proclaim{Theorem}
For each $\zeta\in \C$ the function
\begin{equation}\label{kernel}
K_E(\zeta, z) = \frac i 2
   \frac {E(z)\overline{E(\zeta)}- E^*(z)\overline{E^*(\zeta)}}
            {\pi(z-\bar{\zeta})},
\end{equation}
considered as a function of $z$ {\rm (}extended by continuity to
$z=\bar\zeta${\rm ),} belongs to $H(E)$
and it is the reproducing kernel of $H(E)${\rm ,} i.e.
$$
\langle f, K_E(\zeta,\cdot) \rangle_E=\int_{-\infty}^\infty
\frac{f(t)\overline{K_E(\zeta,t)}}{|E(t)|^2} \,dt =f(\zeta),
$$
for each $f\in H(E)${\rm .}
\endproclaim

For $x\in \R$, write $E(x)$ as
$$
E(x)=|E(x)|e^{-i\varphi(x)},
$$
where $\varphi(x)$ is some continuous function in $\R$ such that
$E(x)e^{i\varphi(x)}$ is real for all $x\in\R$.
If $E(x)\neq 0$, then \eqref{kernel} yields
\begin{equation}\label{major}
\|K_E(x,\cdot)\|^2_E=K_E(x,x)=\frac 1\pi \varphi'(x) |E(x)|^2.
\end{equation}
We say that $\varphi$ is a {\it phase function} of $E$. If there
is a need to distinguish $\varphi$ from phase functions of other
functions in $\overline{\rm HB}$, then we set $\varphi=\varphi_E$.

We shall make use of the following remarkable extension of the Plancherel
identity \cite[p. 55]{dB68}.

\specialnumber{E} \proclaim{Theorem} Let $H(E)$ be a de Branges space and $\varphi$ a phase
function associated with $E${\rm .} Suppose $\alpha$ is a real number
and let $\Gamma=\{\gamma_k\}$ be the {\rm (}possibly finite or one-sided{\rm )} sequence of
real numbers such
that $\varphi(\gamma_k)=\alpha+k\pi${\rm ,} $k\in\Z${\rm .} Then if
$e^{i\alpha}E-e^{-i\alpha}E^*\not\in H(E)${\rm ,} then the normalized
reproducing kernels $K_E(\gamma_k,z)/\|K_E(\gamma_k,\cdot)\|_E$
constitute an orthonormal basis for $H(E)${\rm .} In particular{\rm ,}
$$
\|f\|^2_E=
\sum_k\frac{\pi |f(\gamma_k)|^2}{\varphi'(\gamma_k)|E(\gamma_k)|^2}
$$
for all $f\in H(E)${\rm ;} $e^{i\alpha}E-e^{-i\alpha}E^*\in H(E)$
holds for at most one $\alpha$, modulo $\pi${\rm .}
\endproclaim

The following theorem  will also be crucial \cite{dB60}:
\specialnumber{F} \proclaim{Theorem}  Let $\mu$ be a nonnegative measure on $\R$ and $E$ some
function belonging to $\overline{\rm HB}$ such that $|E(z)|>|E^*(z)|$ for $\Im z>0$
and $E$ has no real zeros{\rm .} Then
$$
\int_{\R} |f(t)/E(t)|^2\, d\mu(t)=\int_{\R}|f(t)/E(t)|^2\, dt
$$
for all $f\in H(E)$
if and only if there exists a bounded holomorphic function $A$ in the
upper half-plane $\C^+$ such that
$ \|A\|_\infty=\sup_{z\in \C^+} |A(z)|\leq 1$ and
$$
\frac{y}\pi\int_{\R}\frac{d\mu(t)}{(t-x)^2+y^2}=\Re\frac{E+E^*A}{E-E^*A}.
$$
\endproclaim

Theorems E and F are due to de Branges for functions $E$ of this kind
(the Hermite-Biehler class: $|E(z)|>|E^*(z)|$ for $\Im z>0$; this excludes
$E^*=cE$ with $|c|=1$, for which $H(E)=\{0\}$; Theorem D holds for all $E\in\overline{\rm HB}$). The functions $E$ to which we
apply them have $H(E)={\rm PW}\neq\{0\}$, and are therefore of this kind
(see Section 3).

We write $H(E)=H(F)$ if $H(E)$ and $H(F)$ coincide when considered as sets
and the norm of $H(F)$ is equivalent to the norm of $H(E)$.
The functions $E\in\overline{\rm HB}$ with the property that ${\rm PW}=H(E)$ have been
characterized in \cite{LS99}.\pagebreak

\section{Main theorem}
\advance\eqcount by 2

If $\Lambda$ is a sampling sequence, then there exists a separated
subsequence $\Lambda'$ which is also sampling
(cf.\ Lemma 3.11 of \cite{Se95}). Hence without loss of generality, we may
restrict our attention to separated sequences $\Lambda$.

\specialnumber{1}\proclaim{Theorem}
A separated sequence $\Lambda$ of real numbers is sampling for ${\rm PW}$ if and
only if there exist two entire functions $E,F\in\overline{\rm HB}$ such that
\begin{itemize}
\item[{\rm (i)}] $H(E)={\rm PW},$
\item[{\rm (ii)}] $\Lambda$ constitutes the zero sequence of $EF+E^*F^*$.
\end{itemize}
\endproclaim 

The following notation will be used repeatedly below: We write
$f\lesssim g$ if there is a constant $K$ such that $f\leq K g$;
we write $f\simeq g$ if both $f\lesssim g$ and $g\lesssim f$.

\demo{Proof} Let us assume first that $\Lambda$ is a
sampling sequence. This means that the Paley-Wiener
space equipped with the norm $\sqrt{\sum_{k} |f(\lambda_k)|^2}$ is a de
Branges space (it is complete since the norm is equivalent to the standard one;
(H2) holds for the same reason; (H1) and (H3) hold because the operations
involved preserve ${\rm PW}$ and, for real $x$, $|(x-\bar\zeta)/(x-\zeta)|=1$ and
$|f^*(x)|=|f(x)|$). Therefore, Theorem C provides us with a function
$E\in \overline{\rm HB}$ such that $H(E)={\rm PW}$ and
$$
\sum_{k} |f(\lambda_k)|^2 =\int_{\R} \frac {|f(t)|^2}{|E(t)|^2}\, dt.
$$
Here $E$ has no real zeros: if $E(x_0)=0$ and $f\in{\rm PW}$ with
$f(x_0)\ne0$, then $\int|f(x+iy)/E(x+iy)|^2dx\gtrsim 1/y$ as $y\to0^+$, so
$f/E\notin H^2$, contradicting $f\in H(E)$. In particular $E(\lambda_k)\neq0$.
Also $|E(z)|>|E^*(z)|$ for $\Im z>0$, since otherwise $E^*=cE$ with $|c|=1$
and $H(E)=\{0\}$. Thus $E$ is in the class to which Theorem~F applies.
Since the measure $\mu=\sum_k |E(\lambda_k)|^2 \delta_{\lambda_k}$ meets the
hypothesis of Theorem F, there exists a bounded holomorphic function
$A$ in the upper half-plane with norm $\|A\|_{\infty}\le 1$ and a real number
$a$  such that
\begin{equation}\label{star}
\frac i\pi \sum_{k} |E(\lambda_k)|^2 \left(\frac1{z-\lambda_k}+
\frac{\lambda_k}{1+\lambda_k^2}\right) +ia=\frac{E(z)+E^*(z)A(z)}{E(z)-E^*(z)A(z)}.
\end{equation}
The series converges, because $\sum_k|E(\lambda_k)|^2/(1+\lambda_k^2)$ is $\pi$
times the Poisson integral of $\mu$ at $z=i$, which is finite by Theorem F.
Note that the right-hand side is a holomorphic function defined in
the upper half-plane, but the left-hand side is a meromorphic
function defined in the whole plane. We denote this meromorphic function
by $M$. The following relationship holds when $\Im z >0$:
$$
A=\frac{M-1}{M+1}\frac{E}{E^*}.
$$
The function $M-1$ has poles at the points $\lambda_k$ as we can
see from \eqref{star}, and it vanishes whenever $E^*$ vanishes (the zeros of
$E^*$ lie in $\Im z>0$, since $E$ has no real zeros, and there
$M-1=2E^*A/(E-E^*A)$ with $E-E^*A\neq0$, because $|E^*A|\le|E^*|<|E|$).
Setting
$$
G(z)= z^{\nu}\prod_{\lambda_k\ne0}\left(1-\frac z{\lambda_k} \right) e^{z/\lambda_k},
$$
where $\nu=1$ if $0\in\Lambda$ and $\nu=0$ otherwise,
we may therefore write
$$
M-1=-\frac{E^*F^*}{G},
$$
with $F$ an entire function.
We also have $M^*=-M$ and $G^*=G$, and consequently $M+1=EF/G$. It
follows that $F^*/F=-A$ in the upper half-plane; moreover
$F=2G/(E-E^*A)$ there, so that $F$ has no
zeros in $\Im z>0$. Since $|A|$ is bounded by 1, it follows that
$F\in\overline{\rm HB}$.

We will now see that $G=(EF+E^*F^*)/2$, which means  that
$\Lambda$ is the zero sequence of $EF+E^*F^*$. We know that
$G=-MG+EF$. If $x\in\R$, then $G(x)$ is real and $M(x)G(x)$ is an
imaginary number. Thus $G$ is the real part of $EF$ for real $z$,
whence $G(z)=(EF+E^*F^*)/2$ for all $z\in \C$.

We next prove the converse implication. The function $E$ has no
zeros on the real axis since ${\rm PW}=H(E)$. We may assume that $F$ has
no real zeros. If this were not the case, then $F=S\widetilde F$,
where $\widetilde F\in\overline{\rm HB}$ without real zeros and $S$ is
an entire function, with real zero set $Z(S)$, and moreover
$S^*=S$ ($S$ is a Weierstrass product over the real zeros of $F$, counted with
multiplicity and with factors of variable genus, since the multiplicities may be
unbounded). Therefore $\Lambda=\Lambda_1\cup Z(S)$ and $\Lambda_1$
will be the zero set of $E\widetilde F+E^*\widetilde F^*$. If we
prove that $\Lambda_1$ is a sampling sequence, then $\Lambda$ is a
also a sampling sequence (the lower inequality passes to the larger set, and the
upper one holds because $\Lambda$ is separated, by the Plancherel-P\'{o}lya
inequality).
\vglue2pt

Given $\alpha\in (0,\pi]$
we let $\Lambda_{\alpha}$ be the sequence of points
$\lambda_{\alpha,k}$ such that
$\varphi_{EF}(\lambda_{\alpha,k})=\alpha + k\pi$, $k\in \Z$.
(The function $\varphi_{EF}=\varphi_E+\varphi_F$ is strictly increasing, since
$\varphi_E'>0$ and $\varphi_F'\ge0$, and it maps $\R$ onto $\R$, because
$\varphi_E(x)\to\pm\infty$ when $x\to\pm\infty$: otherwise the level sets of
$\varphi_E$ would be finite on one side, and Theorem E (for any $\alpha$ but one) would give
$\|g(\cdot-R)\|_E\to0$ as $R\to\infty$ for any fixed $g\in{\rm PW}$, by the
Plancherel-P\'{o}lya inequality (the weights $\pi/(\varphi_E'|E|^2)$ are bounded, by the
estimate of $\varphi_E'|E|^2$ below, and the level sets of $\varphi_E$ are separated, see the
proof of Corollary 1), contradicting $H(E)={\rm PW}$.)
Observe that since $\Lambda$ is the zero sequence of $EF+E^*F^*$, and $E$ and
$F$ have no real zeros,
we have $\Lambda=\Lambda_{\pi/2}$. For $\alpha\ne \pi/2$ the
sequence $\Lambda_\alpha$ is interlaced with the sequence
$\Lambda$. Hence since $\Lambda$ is a separated sequence,
$\Lambda_\alpha$ can be expressed as
the union of two separated sequences, each
with separation constant not smaller than that of $\Lambda$. Thus
the Plancherel-P\'{o}lya inequality \cite[pp. 96--98]{Yo80}
implies that
\begin{equation}\label{planpol}
\sum_k |f(\lambda_{\alpha,k})|^2\leq C \|f\|^2_{{\rm PW}},
\end{equation}
with $C$ independent of $\alpha$. Moreover, for all values of
$\alpha$ except possibly one, Theorem E implies that for every
$g\in H(EF)$ the following relation holds:
$$
\int_\R \frac {|g(t)|^2}{|E(t)F(t)|^2}\, dt=\sum_{k}
\frac{\pi|g(\lambda_{\alpha,k})|^2} {|E(\lambda_{\alpha,k})|^2
|F(\lambda_{\alpha,k})|^2 \varphi'_{EF} (\lambda_{\alpha,k})}.
$$
For every $f\in H(E)$ the function $g=Ff$ belongs to $H(EF)$ with
$\|f\|_{H(E)}=\|Ff\|_{H(EF)}$. Therefore, for every $f\in {\rm PW}$ we have
\vglue-9pt
\begin{equation}\label{equiv}
\|f\|^2_{{\rm PW}}\simeq \int_{\R} \Bigl|\frac{f(t)}{E(t)}\Bigr|^2\, dt
= \int_{\R}\Bigl|\frac {g(t)} {E(t)F(t)}\Bigr|^2\, dt=\sum_{k}
\frac{\pi|f(\lambda_{\alpha,k})|^2}
{|E(\lambda_{\alpha,k})|^2\varphi'_{EF}(\lambda_{\alpha,k})}.\qquad
\end{equation}
Since $H(E)={\rm PW}$, $E$ does not vanish for $x\in \R$ and therefore
by \eqref{major}
$$
1=\sup_{{f\in {\rm PW}\atop \|f\|^2_{{\rm PW}}\le 1}}|f(x)|^2\simeq
\sup_{{f\in H(E)\atop \|f\|^2_E\le 1}} |f(x)|^2=K_E(x,x)=\frac
1\pi \varphi_E'(x)|E(x)|^2
$$
for $x\in\R$. Since
$\varphi'_{EF}=\varphi'_E+\varphi'_F\ge\varphi'_E$ (here $\varphi_F'\ge0$,
because $\log|F^*|-\log|F|\le0$ in $\Im z>0$, with equality on $\R$), we obtain from
\eqref{equiv} the inequality
$$
\|f\|^2_{{\rm PW}}\leq C \sum_{k}|f(\lambda_{\alpha,k})|^2,
$$
where the constant $C$ does not depend on $\alpha$.
This inequality may fail for one $\alpha\in (0,\pi]$. We may assume it fails
for $\alpha=\pi/2$, because otherwise we have proved that $\Lambda$ is
sampling. Take a sequence $\alpha_n\to \pi/2$. Then
$\lambda_{\alpha_n,k}\to\lambda_k$ for all $k$. Let $f\in{\rm PW}$. By the
local form of the Plancherel-P\'{o}lya inequality (the discs of radius $q/2$
about the points of $\Lambda_\alpha$, $q$ the separation constant of $\Lambda$,
overlap at most twice, and $|f|^2$ is subharmonic),
$$\sum_{|\lambda_{\alpha,k}|>R}|f(\lambda_{\alpha,k})|^2\lesssim
\int\!\!\int_{|x|>R-q,\,|y|<q}|f(x+iy)|^2\,dx\,dy$$
uniformly in $\alpha$, and the right-hand side tends to $0$ when $R\to\infty$.
Hence $\sum_k|f(\lambda_{\alpha_n,k})|^2\to\sum_k|f(\lambda_k)|^2$, and the
inequality holds also for $\alpha=\pi/2$.
\enddemo

If $\Lambda$ is a complete interpolating sequence, then there
exists an $E\in \overline{\rm HB}$ such that $H(E)={\rm PW}$ and $\Lambda$
constitutes the zero sequence of $E+E^*$. This follows from
Theorems C, D and E: a complete interpolating sequence is sampling (by the open
mapping theorem), so that Theorem C applies; the kernels $K_E(\lambda_k,\cdot)$
are then orthogonal and complete in $H(E)$, so that by Theorem D the numbers
$\varphi_E(\lambda_k)$ differ by multiples of $\pi$; after replacing $E$ by
$e^{i\theta}E$ we get $\Lambda=\{\varphi_E\in\pi/2+\pi\Z\}$, the zero set of
$E+E^*$. Therefore, we may view the function $F$ of
Theorem 1 as accounting for the ``redundancy'' in $\Lambda$. There
is a trivial case in which this is particularly transparent:
Suppose $\Lambda=\Lambda'\cup (\Lambda\setminus\Lambda')$, where
$\Lambda'$ is a complete interpolating sequence. (Such a
decomposition holds whenever $D^-(\Lambda)>1$; cf.\ Theorem~2.3 of
\cite{Se95}.) Then the condition of Theorem 1 is met if we
choose $F$ to be the generating function of the ``redundant''
sequence $\Lambda\setminus \Lambda'$, i.e., if we set
$$
F(z)= \prod_{\lambda_k\in\Lambda\setminus\Lambda'}
\left(1-\frac z{\lambda_k} \right) e^{z/\lambda_k}
$$
(with the factor $z$ in place of $(1-z/\lambda_k)e^{z/\lambda_k}$ when $\lambda_k=0$).

A somewhat different way of seeing $F$ as measuring the
``redundancy" in $\Lambda$ is as follows. The notion of a complete
interpolating sequence may be extended in a natural way: We say
that $\Lambda$ is a complete interpolating sequence for a de
Branges space $H$ if the interpolation problem $f(\lambda_k)=a_k$
has a unique solution $f\in H$ whenever
$$
\sum_k \frac{|a_k|^2}{K(\lambda_k,\lambda_k)}<\infty.
$$
Then Theorem 1 says that if $\Lambda$ is sampling for {\rm PW}, then
$\Lambda$ is a complete interpolating sequence for $H(EF)$, and
$H(E)$ (with $H(E)={\rm PW}$) is isometrically embedded into $H(EF)$ by
the map $f\to F f$. This interpretation of Theorem 1 has an
interesting relation to the result of \cite{Se95} which says
that we cannot in general obtain from $\Lambda$ a complete
interpolating sequence for ${\rm PW}$ by making {\it the sequence
thinner}: Instead, we can make {\it the space bigger} so that
$\Lambda$ becomes a complete interpolating sequence for the bigger
space.

We note in passing that Theorem 1 solves the following problem of
norm equivalence for ${\rm PW}$: Which nonnegative measures $\mu$ on
$\R$ have the property that the $L^2(\R,d\mu )$-norm yields a norm
for ${\rm PW}$ equivalent to the standard $L^2(\R)$-norm? To apply
Theorem~1 to this problem, we need a way of associating sampling
sequences with such measures. Given a nonnegative measure $\mu$
and two positive numbers $r,\delta$, define
$$
\Lambda_\mu(\delta,r)=\{kr:\ k\in\Z \ \hbox{and}\
\mu([kr,(k+1)r))\geq \delta\}.
$$
Then the Bernstein and the Plancherel-P\'{o}lya inequalities imply
(see \cite{Or98}):
\nonumproclaim{Proposition}
The $L^2(\R,\mu)$\/{\rm -}\/norm yields an equivalent norm for ${\rm PW}$ if and
only if the following holds\/{\rm :}
\begin{itemize}
\item[{\rm (i)}]There
exists a positive constant $C$ such that $\mu([x,x+1))\leq C$ for
all $x\in\R${\rm .}
\item[{\rm (ii)}]  For all sufficiently small $r>0$ there exists a
$\delta=\delta(r)>0$ such that $\Lambda_\mu(r,\delta)$ is sampling
for ${\rm PW}${\rm .}
\end{itemize}
\endproclaim

Problems about norm equivalence for spaces of entire functions of
exponential type in one and several variables have been studied by
many authors. See, e.g., \cite{Li65} and \cite{LLS92} for an
extensive historical account. 

\section{Landau's inequality and the John-Nirenberg theorem}

The following corollary says that a sampling sequence is ``everywhere
denser'' than some complete interpolating sequence.

\specialnumber{1}\proclaim{{C}orollary}\label{cis}
If $\Lambda$ is a separated sampling sequence for ${\rm PW}${\rm ,} then there
exists a complete interpolating sequence
$\Gamma=\{\gamma_k\}_{k=-\infty}^\infty$ such that for every $k\in
\Z$ there is at least one point $\lambda\in \Lambda$ such that
$\gamma_k\leq \lambda<\gamma_{k+1}${\rm .}
\endproclaim

\demo{Proof}
If $\Lambda$ is a sampling sequence, let $E$ and $F$ be as in Theorem 1 and
write $F=S\widetilde F$ as in its proof ($S=1$ if $F$ has no real zeros). Then $\Lambda$
contains the zero set $\Lambda_1$ of $E\widetilde F+E^*\widetilde F^*$, which consists of
those points such that $\varphi_{E\widetilde F}(\lambda)=\pi/2+k\pi$ for some
$k\in \Z$. On the other hand, $\varphi_{E}$ is an increasing
function which grows more slowly than $\varphi_{E\widetilde F}$: the difference
$\varphi_{\widetilde F}$ is nondecreasing, so that $\varphi_{E\widetilde F}$ increases by at
least $\pi$ between consecutive points of the sequence $\Gamma_\alpha$ below. Thus the
sequence $\Gamma_\alpha$ which consists of those points such that
$\varphi_E(\gamma)=\alpha+k\pi$ for some $k\in \Z$, has the
claimed property. Since $H(E)={\rm PW}$, Theorem E implies that
$\Gamma_\alpha$ is a complete interpolating sequence for all
values of $\alpha$ except possibly one, mod $\pi$. (The sequences $\Gamma_\alpha$
are separated, as required in Theorem A: the function $f_\gamma\in{\rm PW}$ with
$f_\gamma(\gamma)=1$ and $f_\gamma=0$ on $\Gamma_\alpha\setminus\{\gamma\}$
has $\|f_\gamma\|^2_{\rm PW}\lesssim\|f_\gamma\|_E^2=\pi/(\varphi_E'(\gamma)|E(\gamma)|^2)\lesssim1$,
and Bernstein's inequality gives $|f_\gamma'|\lesssim\|f_\gamma\|_{\rm PW}$.)
Since the
choice of $\alpha$ is at our disposal, the proof is complete.
\enddemo

Let us return to Landau's inequality, which is stated as Theorem B
above. To begin with, let us show by an example that the
inequality is best possible. Set
$\Lambda=\{k+\log^+|k|\}_{k\in\Z}$. This sequence is sampling,
because the function
$$ F(z)=z\lim_{R\to\infty}\prod_{0<|\lambda_k|<R}
\left(1-z/\lambda_k\right) $$ is a sine-type function; i.e.,
$|F(z)|e^{-\pi |\Im z|}\simeq 1$ for $|\Im z|>1$. (Indeed,
$\log|F(z)|-\pi|\Im z|$ is the potential of the measure $d(n_\Lambda(t)-t)$, and
$n_\Lambda(t)-t=-\log^+|t|+O(1)$, where the bounded term is a sawtooth with unit jumps at
the points $\lambda_k$. Boundedness alone does not suffice, but for $|\Im z|\ge1$ the point
masses at the $\lambda_k$ may be replaced, as in Lemma 3 of \cite{OS99}, by a density $m\simeq1$
having the same mass and barycentre on each interval between the midpoints of consecutive
points $\lambda_k$; then $\int_0^t(m-1)=-\log^+|t|+c+O((1+|t|)^{-1/2})$. Since the measure
$d\log^+|t|$ is odd, only the conjugate function of $\log|t|$ comes into play, and this is
bounded; the remainder contributes $O(1)$ for $|\Im z|\ge1$.)
Thus a classical
theorem of Levin \cite[p. 250]{HNP81} implies that $\Lambda$ is
a complete interpolating sequence, and hence sampling. (Alternatively, note that this
means $v=0$ in Theorem A.) On the other hand, it is clear that
$n_\Lambda(R)=R-\log R+O(1)$ and $-n_\Lambda(-R)=R+\log R+O(1)$
when $R\to\infty$.

\nonumproclaim{Claim} Landau\/{\rm '}\/s inequality is a consequence of
the John\/{\rm -}\/Nirenberg theorem for $\BMO$ functions{\rm .}
\endproclaim
\advance\eqcount by 5

\demo{Proof}
Suppose $\Lambda$ is a sampling sequence, and let $\Gamma$ be an associated
complete interpolating sequence as described in Corollary 1. By Corollary 1
and Theorem A,
$$
n_{\Lambda}(b)-n_\Lambda(a)\ge n_\Gamma(b)-n_\Gamma(a)-1\ge b-a +
h(b)-h(a)-1,
$$
where $h\in {\rm BMO}$ of a special form. (For this proof we will only
need the fact that $h\in {\rm BMO}$.) Set $I=[a,b]$. Then the triangle
inequality in the form
$$
|h(b)-h(a)|\le |h(b)-h_I| + |h(a)-h_I|,
$$
gives
\begin{equation}\label{ineq}
n_{\Lambda}(b)-n_\Lambda(a)\ge b-a -2\max_{t\in I}|h(t)-h_I|-1.
\end{equation}
We may assume that $|I|>2$, since $h=n_\Gamma-t$ has variation $\lesssim1$ on
intervals of length $\le2$, and such intervals only change the constant $B$ below.
Define
$$
J=\{t\in I; |h(t)-h_I|\ge A \log |I|\}.
$$
Then the John-Nirenberg theorem (cf.\ \cite{JN61}, \cite[p.
230]{Ga81}) implies
$$
|J|\le C e^{-\frac {cA\log
|I|}{\|h\|_*}} |I| \le C/|I|
$$
for some sufficiently big $A$, with $C,c$ absolute positive
constants. But $h'(t)$ is bounded, in fact equal to $-1$ for
$t\not\in\Gamma$, and so if the right-hand side is smaller than
the separation distance of $\Gamma$,  then every $t\in I$ is within that
distance of some point of $I\setminus J$, and consequently $|h(t)-h_I|\leq
A\log|I|+B$ for all $t\in I$. Plugging this into \eqref{ineq}, we
obtain Landau's inequality.
\enddemo

It should be stressed that we do not pretend to have given an
easier proof of Landaus's inequality, as the John-Nirenberg
theorem is admittedly a more sophisticated result. However, we
find the link between the two results quite intriguing.

\section{Sampling sequences and approximation\\ of subharmonic functions}
\advance\eqcount by 6

In Section 3, we explained how Theorem 1 could be
interpreted in the trivial case in which there exists a
subsequence $\Lambda'\subset \Lambda$ which is a complete
interpolating sequence. The purpose of this section is to
demonstrate the applicability of our condition when
$D^-(\Lambda)=1$ and no such subsequence $\Lambda'$ exists.

To obtain positive results from Theorem 1, one needs information
about those $E\in\overline{\rm HB}$ for which $H(E)={\rm PW}$. As mentioned
at the end of Section 2, such functions are described in
\cite{LS99}. However, we will not use this description, which
involves a rather delicate statement about the location of the
zeros of $E$. In fact, we will obtain precise
results by using only the fact that $H(E)={\rm PW}$ if
$|E(z)|\simeq e^{\pi \Im z}$ for $\Im z\geq 0$. (If $f\in H(E)$, then
$|f(z)|\lesssim\|f\|_E|E(z)|/\sqrt{\Im z}$ for $\Im z>0$ and $f\in L^2(\R)$ by
Fatou's theorem, while a Phragmén-Lindel\"of argument in the strip $|\Im z|\le1$
shows that $f$ has exponential type at most $\pi$; conversely, if $f\in{\rm PW}$,
then $\int|f(x+iy)/E(x+iy)|^2dx\lesssim\int|f(x+iy)|^2e^{-2\pi y}dx\le\|f\|^2_{\rm PW}$.)

We will assume $\psi\in {\mathcal C}^1(\R)$ is a nondecreasing
function satisfying\break $\psi(\infty)-\psi(-\infty)=\infty$ and
$\psi'(x)=o(1)$ as $|x|\to\infty$. We associate with $\psi$ a
sequence $\Lambda(\psi)=\{\lambda_k\}_{k\in\Z}$ defined by
$\lambda_k+\psi(\lambda_k)=k$. In other words, with $[x]$ denoting
the integer part of the real number $x$ and $\psi(0)=0$, this
means that $n_{\Lambda(\psi)}(t)=[t+\psi(t)]$. All sequences
$\Lambda(\psi)$ are sets of uniqueness of infinite excess, but
none of them contain subsequences which are complete interpolating
sequences; cf.\ Theorem 2.7 in \cite{Se95}. We introduce the
potential
$$
U_\psi(z)=\int_{-\infty}^\infty [\log|1-z/t|+\chi(t)\Re z/t]\, d\psi(t),
$$
where $\chi$ is the characteristic function of $\{|t|>1\}$; the integral
converges absolutely, since $\psi'$ is bounded. (The cut-off $\chi$ replaces the
principal value sense, which need not exist at $t=0$; the two potentials differ
by $c\Re z$ with $c$ real, when the principal value exists, and this is
immaterial below, because $f\mapsto fe^{-cz}$ preserves $\overline{\rm HB}$ and the
phase function.) The same convention is used for all potentials below.
This function is
subharmonic because $\psi'(t)\geq 0$.

We will now demonstrate how our sampling condition is related to
the problem of approximating $U_\psi$ by $\log |f|$, with $f$ an
entire function.

\specialnumber{2}\proclaim{{C}orollary} A sequence $\Lambda(\psi)$ is sampling for ${\rm PW}$ if there
exists an $f\in\overline{\rm HB}$ such that $\varphi_f'(x)=o(1)$ when
$|x|\to\infty$ and
\begin{equation}\label{appr}
|U_\psi(z)-\log|f(z)|| \lesssim 1  \ \ \hbox{for}  \ \ \Im z\geq
0.
\end{equation}
\endproclaim

\demo{Proof}
If we could find $e\in \overline{\rm HB}$ such that
$$
\varphi_e(x)=\pi x+\pi \psi(x)-\varphi_f(x),
$$
we would be done, because then $\Lambda$ would constitute the zero
sequence of $ef+e^*f^*$, and $|e(z)|\simeq e^{\pi \Im z}$ for $\Im
z\geq 0$. In general, however, it is impossible to find such an
$e$. Instead, we will appeal to the following perturbation
argument: If a separated sequence $\Gamma=\{\gamma_k\}$ is
sampling, then $\Gamma'=\{\gamma_k+\delta_k\}$ is sampling
whenever all $\gamma_k+\delta_k$ are distinct and $\delta_k\to 0$ as $|k|\to
\infty$. (This follows from Lemma III in \cite{DS52} and the
elementary fact that $\Gamma$ is sampling if $\delta_k=0$ except
for finitely many integers $k$; here finitely many points are moved to new
positions, distinct from the others, whereas removing points is not allowed.
The elementary fact is proved by a normal family argument.) Thus it is enough to construct an
$E\in \overline{\rm HB}$ such that
\begin{equation}
\varphi_E(x)-\pi x-\pi \psi(x)+\varphi_f(x)=o(1) \label{enough}
\end{equation}
when $|x|\to\infty$ and $|E(z)|\simeq e^{\pi \Im z}$ for $\Im
z\geq 0$, because then the perturbation argument applies with the
zero sequence of $Ef+E^*f^*$ playing the role of $\Gamma$ and
$\Gamma'=\Lambda$. (This zero sequence is separated, being a level set of
$\varphi_{Ef}=\varphi_E+\varphi_f$ and $\varphi_{Ef}'=\varphi_E'+\varphi_f'\lesssim1$,
since $\varphi_E'\simeq|E|^{-2}\simeq1$ on $\R$ and $\varphi_f'$ is bounded.)
\vglue2pt

By \eqref{appr}, $f$ has no zeros in $\Im z\ge0$, so that $\log f$ is analytic in
$\Im z>0$ and continuous up to $\R$, and $\Im\log f(x)=-\varphi_f(x)+{\rm const}$ on
$\R$. For a bounded function $g$ set
$U_g(z)=\int[\log|1-z/t|+\chi(t)\Re z/t]\,g(t)\,dt$. By the Cauchy-Riemann
equations, $\log|f|$ and $U_{\varphi_f'}/\pi$ have the same normal derivative
$\varphi_f'$ on $\R$, so their difference extends by reflection to an entire
harmonic function, which is of growth $o(|z|^2)$ by \eqref{appr}, hence of the
form $b\Re z+c$ with $b,c$ real:
$$
\log|f(z)|=\frac1\pi U_{\varphi_f'}(z)+b\Re z+c,\qquad \Im z\ge0 .
$$
Put $\omega(x)=x+\psi(x)-\varphi_f(x)/\pi$ and $U_\omega=U_{\omega'}$. Since
$U_1(z)=\pi\Im z$ for $\Im z\ge0$, we get
$U_\omega=\pi\Im z+U_\psi-\log|f|+b\Re z+c=\pi\Im z+b\Re z+O(1)$ by \eqref{appr}.
Now $\omega'=1+\psi'-\varphi_f'/\pi\to1$ at infinity, because $\varphi_f'=o(1)$.
We may assume that $\omega'(x)\simeq 1$ on $\R$ and that $b=0$: add to
$\omega'$ a compactly supported continuous function $h$ such that $\omega'+h\simeq1$,
$\int h=0$ and $\int_{|t|>1}h(t)\,dt/t=-b$ (this is possible by adding a bump far
out, where $\omega'\simeq1$). The condition $\int h=0$ keeps the phase unchanged at
$\pm\infty$, so that \eqref{enough} is not affected, and the second condition cancels
the linear term $\Re z\int_{|t|>1}h(t)\,dt/t$ of the potential of $h$.
Thus we may assume $U_\omega(z)=\pi\Im z+O(1)$ for $\Im z\ge0$. Partition the real line into a
sequence of disjoint intervals $I_k=[x_k,x_{k+1}]$, $k\in\Z$, with
$x_0=0$,  such that
$$
\int_{I_k}\omega'(t)dt=1
$$
for all $k$, and choose $\gamma_k\in I_k$ so that
$$
\gamma_k=\int_{I_k} t \omega'(t) dt.
$$
We set
$$
A(z)=\lim_{R\to\infty}e^{c_Rz}\prod_{|\gamma_k|<R}(1-z/\gamma_k)
$$
for suitable real constants $c_R$ (the symmetric products need not converge by
themselves), and $\Gamma=\{\gamma_k\}$. Then by Lemma 3 in \cite{OS99},
$$
|A(z)|e^{-U_\omega(z)}\simeq \min(1,\dist(z,\Gamma)),
$$
so that $|A(z)|\simeq e^{\pi\Im z}\min(1,\dist(z,\Gamma))$ for $\Im z\ge0$.
Now choose two polynomials $P$ and $Q$ of the same degree and with
only real zeros, with $P/Q\to1$ at infinity, such that
$B(z)=A(z-1/2)P(z)/Q(z)$ is an entire
function and the zeros $\beta_k$ of $B$ and the zeros $\gamma_k$ of $A$ are interlaced
(the gaps of $\Gamma$ tend to $1$, so that $Q$ has to cancel, and $P$ to
replace, only finitely many zeros of $A(z-1/2)$). We set $E=A-iB$.
Then $|E|^2+|E^*|^2=2(|A|^2+|B|^2)$ and $|E^*|\le|E|$ in $\Im z\ge0$ (see below), so that
$|A|^2+|B|^2\le|E|^2\le2(|A|^2+|B|^2)$ there, and since
$|B(z)|\simeq|A(z-1/2)|\simeq e^{\pi\Im z}\min(1,\dist(z,\Gamma+1/2))$ for large $|z|$, while
$\Gamma$ and $\Gamma+1/2$ are at distance $\simeq1/2$ from each other, we obtain
$|E(z)|\simeq e^{\pi\Im z}$ for $\Im z\ge0$.

We claim that $E\in\overline{\rm HB}$. (Interlacing zeros and exponential type
are not sufficient for this by themselves: $\cos\pi z$ and $e^{az}\sin\pi z$, $a\ne0$ real,
have interlacing zeros, but $|\cos\pi z\mp ie^{az}\sin\pi z|$ is larger at $iy$ than at $-iy$
for some $y>0$ and smaller for others. So we argue directly.)
Let $A_R(z)=\prod_{|\gamma_k|<R}(1-z/\gamma_k)$ and let $B_R=A_R(\cdot-1/2)P/Q$,
for $R$ large. These are real polynomials, and
$B_R/A_R=\prod_{|k|<R}(z-\beta_k)/(z-\gamma_k)$, with interlacing zeros
$\gamma_{-R}<\beta_{-R}<\gamma_{-R+1}<\dots$. If $\Im z>0$, the argument of
each factor $(z-\beta_k)/(z-\gamma_k)$ lies in $(0,\pi)$, and by interlacing
the sum of these arguments telescopes to a number in $(0,\pi)$; hence
$\Im(B_R/A_R)>0$ in $\Im z>0$. Therefore $|A_R+i\lambda B_R|<|A_R-i\lambda B_R|$ in
$\Im z>0$ for every $\lambda>0$, that is, $A_R-i\lambda B_R\in\overline{\rm HB}$. Now
$e^{c_Rz}(A_R-ie^{-c_R/2}B_R)\to A-iB=E$ locally uniformly, and $\overline{\rm HB}$ is
closed under locally uniform limits (by Hurwitz' theorem), so $E\in\overline{\rm HB}$.
(No choice of sign of $P$ is needed.)

It remains to show  that $E$ satisfies \eqref{enough}. Consider
the sequence of functions $A_k(z)=A(z+x_{2k})$. By the estimate for $|A|$ above,
these functions are uniformly bounded by $Ce^{\pi|\Im z|}$, so that they form
a normal family; since $\omega'\to1$, the zeros of $A_k$ tend to $\Z+1/2$, and
any limit function is a real entire function of exponential type at most $\pi$ with
zero set $\Z+1/2$, which is bounded on $\R$. Dividing by $\cos\pi z$ and
applying Liouville's theorem, we find that the limit is $c\cos\pi z$ with
$c$ real, and $|c|\simeq1$ by the lower estimate. Hence there exists a sequence $c_k\simeq 1$ such
that $A_k(x)-c_k \cos\pi x\to 0$ uniformly on compact subsets of
the real line. Similarly, if we set $B_k(z)=B(z+x_{2k})$, we
obtain that $B_k(x)-c_k\sin\pi x\to 0$ uniformly on compact
subsets of the real line. The same holds when $k\to-\infty$.
Then it follows from the construction of
$E$ that
$$
E(x)=|E(x)|e^{-i(\pi(x+\psi(x))-\varphi_f(x)+o(1))}
$$
when $|x|\to\infty$, with $|E(x)|\simeq 1$, after replacing $E$ by $e^{i\theta}E$
if necessary. (The constant is the same at $+\infty$ and at $-\infty$: the zeros
$\gamma_k$ of $A=(E+E^*)/2$ are the points where $\varphi_E\in\pi/2+\pi\Z$, in
increasing order, so that $\varphi_E(\gamma_k)=\pi/2+\pi(k+n_0)$ for all $k$,
with one integer $n_0$.)
\enddemo

To solve the approximation problem for $U_\psi$, we shall adapt an
approach from a recent paper by Lyubarskii and Malinnikova
[LM01], which contains rather conclusive results about
the possibility of approximating an arbitrary subharmonic function
by $\log|f|$ with $f$ an entire function.

To ease the exposition, we set $\psi(x)=0$ for $x\leq 0$. (It is
not difficult to modify the construction below to get a
corresponding result when both $\psi(\infty)=\infty$ and
$\psi(-\infty)=-\infty$.) Let $\{t_n\}_{n=0}^\infty$ be the
sequence such that $t_0=0$ and $\psi(t_n)=n$, $n=1,2,3,...$; set
$d_n=t_{n}-t_{n-1}$. We will say that $\psi$ induces a
{\it logarithmically regular partition} if $d_n\simeq d_{n+1}$
and
$$
\sup_{x>0}\sum_{x/2<t_n<2x}\frac{d_n^2}{(x-t_n)^2+d_n^2}<\infty.
$$

\specialnumber{2}\proclaim{Theorem}
Suppose $\psi(x)=0$ for $x\leq 0${\rm .} Then
\begin{itemize}
\item[{\rm (i)}] If
$\psi'(x)=1/O(x)$ when $x\to\infty${\rm ,} $\psi'$ is eventually nonincreasing{\rm ,}
and $\psi$ induces a
logarithmically regular partition{\rm ,} then $\Lambda(\psi)$ is
sampling for ${\rm PW}${\rm .}
\item[{\rm (ii)}] If $\psi'(x)=o(1/x)$ when $x\to\infty${\rm ,} then $\Lambda(\psi)$ is not
sampling for ${\rm PW}${\rm .}
\end{itemize}

\endproclaim

We note that $\psi(x)=\sqrt{x}$ for $x\geq 0$ corresponds to the
example of \cite{Se95} showing the existence of a sampling
sequence containing no subsequence which is a complete interpolating
sequence.

(The assumption that $\psi'$ is eventually nonincreasing is a regularity
assumption, which is used only through the bound $\psi'\lesssim1/d_n$ on
$[t_{n-1},t_n]$; it holds, e.g., for $\psi(x)=\sqrt x$, and (i) remains true with this
bound in place of the monotonicity of $\psi'$. Whether (i) holds without an assumption of
this kind is not known. On the other hand, the hypothesis $\psi'(x)=1/O(x)$ is not needed:
the proof only uses $t_n\lesssim t_{n-1}$, which follows from $d_n\simeq d_{n+1}$, and which,
together with the monotonicity of $\psi'$, implies $\psi'(x)=1/O(x)$.)

Case (i) is proved by means of Corollary 2, while case (ii) is
proved by a direct argument. The point of case (ii) is to show
that the growth condition $\psi'(x)=1/O(x)$ is critical. Thus
Theorem~2 indicates that the condition of Corollary 2 is in a
sense close to being necessary.

\demo{Proof}
We consider first statement (i). According to Corollary 2, we need to
show that the assumptions on $\psi$ ensure the existence of a
function $f\in\overline{\rm HB}$ such that \eqref{appr} holds when
$\Im z\geq 0$.

The zeros of $f$ are determined as follows. Define
$r_n\in(t_{n-1},t_n)$ by
\begin{equation} \label{logcent} \log
r_n=\int_{t_{n-1}}^{t_n}\log t \ d\psi(t),
\end{equation}
and then
$$
z_n=r_n e^{-ic d_n/r_n}
$$
with $c>0$ so small that
$cd_n/r_n\leq \pi/4$ for all $n$. We choose $f$ in such a way that it satisfies
$$
\log|f(z)|=\sum_{n=1}^\infty\int_{t_{n-1}}^{t_n}\left ( \log\left |
1-\frac z {z_n} \right |+ \chi(t)\Re \frac zt\right)\, d\psi(t).
$$
Clearly $f$ belongs to $\overline{\rm HB}$.
We set $V=U_\psi-\log|f|$; it satisfies the following relation:
\begin{equation}\label{deff}
V(z)=\sum_{n=1}^\infty\int_{t_{n-1}}^{t_n}\left ( \log\left |
1-\frac z t \right |-
       \log\left | 1-\frac z {z_n} \right |
                         \right )\, d\psi(t)=
\sum_{n=1}^\infty j_n(z).\qquad
\end{equation}
We see that it is enough to prove that the series converges
uniformly on compact sets in $\C$, and that $V(z)=O(1)$ for $\Im z
\geq 0$.

The proof is a simplified version of an argument from
\cite{LM01}. We shall therefore be brief and refer to
\cite{LM01} for further details.

Given $z\in\C$, we let $n(z)$ be a positive integer $n$ such that
$r_{n-1}<|z|\leq r_n$. Fix a positive integer $m$, depending only on the
constant in $d_n\simeq d_{n+1}$, such that the cells with index outside
$[n(z)-m,n(z)+m]$ are at distance $\gtrsim d_n$ from $z$. Then if $\Im z \geq 0$, it is plain that
the bound $\psi'\lesssim1/d_n$ on $[t_{n-1},t_n]$ (which follows from the
assumption that $\psi'$ is eventually nonincreasing: for large $n$,
$\psi'\le\psi'(t_{n-1})\le1/d_{n-1}\lesssim1/d_n$ on $[t_{n-1},t_n]$, because the mean
of $\psi'$ over $[t_{n-2},t_{n-1}]$ is $1/d_{n-1}$) ensures that
$$
\sum_{n=n(z)-m}^{n(z)+m}
j_n(z)\simeq 1.
$$
Next let $n^{-}(z)$ be the positive integer $n$
such that $r_{n-1}<|z|/2\leq r_n$, and $n^{+}(z)$ be the positive
integer $n$ such that $r_{n-1}<2|z|\leq r_n$. Then for $z\in \C$
we see that
$$
\sum_{n=n^+(z)+1}^{\infty} j_n(z)\simeq 1;
$$
here, after the splitting of $j_n$ used for $V_1$ and $V_2$ below, the first order term of
the expansion used for $V_1$ cancels by \eqref{logcent}, and what remains is
$\lesssim|z|d_n/t_{n-1}^2$; since $t_n\lesssim t_{n-1}$, these bounds have sum
$\lesssim|z|\int_{t_{n^+(z)}}^\infty dt/t^2\lesssim1$. (A bound
$|j_n(z)|\lesssim|z|^2/t_n^2$, as for the factors of a canonical product, would not suffice.)
After observing that by \eqref{logcent} we may write
$$
j_n(z)=\int_{t_{n-1}}^{t_n}\left ( \log\left | 1-\frac t z
\right |-
       \log\left | 1-\frac  {z_n} z\right |
                         \right )\, d\psi(t),
$$
we obtain
$$
\sum_{n=1}^{n^{-}(z)-1} j_n(z)\simeq 1
$$
in a similar way (now $|j_n(z)|\lesssim d_n/|z|$, and these bounds have sum $\lesssim1$).
Thus we have in particular established the
uniform convergence on compact sets.

We set
$$
N(z)=\{n^-(z), n^-(z)+1, ..., n^+(z)-1, n^+(z)\}
\setminus\{n(z)-m,...,n(z)+m\}
$$
(the set of ``essential indices'') and split the corresponding sum
into two parts,
\begin{eqnarray*}
V_1(z)+V_2(z)
&=&\sum_{n\in N(z)}\int_{t_{n-1}}^{t_n}\left ( \log\left | 1-\frac z
t \right |-
       \log\left | 1-\frac z {r_n} \right |
                         \right )\, d\psi(t) \\
    &&   +\ \sum_{n\in N(z)}\left ( \log\left |
z-r_n \right |-
       \log\left | z- z_n \right |
                         \right )
\end{eqnarray*}
now assuming that $\Im z\geq 0$.

To estimate $V_1$, introduce the function
$$
L(\omega)=\log(1-ze^{-\omega}),
$$
which for each $n\in N(z)$ is analytic in a domain containing
those $\omega$ such that $e^\omega\in [t_{n-1},t_n]$. Therefore,
$$
L(\omega)-L(\omega_n)=(\omega-\omega_n)L'(\omega_n)
    +\int_{\omega_n}^\omega L''(\sigma)(\omega-\sigma)d\sigma=
    (\omega-\omega_n)L'(\omega_n)+Q_n(z,t),
$$
where $t=e^\omega\in [t_{n-1},t_n]$ and $e^{\omega_n}=r_n$. By
\eqref{logcent}, we obtain
$$
V_1(z)=\Re \sum_{n\in
N(z)}\int_{t_{n-1}}^{t_n} Q_n(z,t)d\psi(t).
$$
A direct estimate gives
$$
\sup_{t\in [t_{n-1},t_n]}|Q_n(z,t)|\lesssim
\frac{d_n^2}{|z-r_n|^2}.
$$
 Using the assumption that $\psi$
induces a logarithmically regular partition, we obtain that
$V_{1}(z)\simeq 1$.

To deal with $V_2$, we write
$$
\log\left | z-r_n \right |-
       \log\left | z- z_n \right |=\Re
       \int_{z_n}^{r_n}\frac{d\zeta}{\zeta-z}.
$$
Integrating along the arc $\zeta=r_n e^{-i\theta}$, $0\leq
\theta\leq cd_n/r_n$, we obtain from this
$$
 \left|\log\left | z-r_n \right |-
       \log\left | z- z_n \right |\right|\lesssim
       \frac{d_n^2}{|z-r_n|^2}+\frac{d_n\Im z}{|z-r_n|^2}.
$$
(The second term is the first order term: along the arc,
$\Re\,d\zeta/(\zeta-z)=\Im\bigl(z/(\zeta-z)\bigr)d\theta$, which is of size $r_n\Im z/|z-r_n|^2\,d\theta$.)
Using again the assumption that $\psi$ induces a logarithmically regular
partition, the sum of the first terms is $\lesssim1$. The sum of the second terms
is a Riemann sum of a Poisson integral: since $|z-t|\simeq|z-r_n|$ for
$t\in[t_{n-1},t_n]$ when $n\in N(z)$,
$\sum_{n\in N(z)}d_n\Im z/|z-r_n|^2\lesssim\int_{\R}\Im z/|z-t|^2\,dt\le\pi$.
Hence $V_{2}(z)\simeq 1$, and we are
done with the proof of \eqref{deff}.

It remains to check that $f$ satisfies the other hypothesis of Corollary 2,
$\varphi_f'(x)=o(1)$. Since $f$ has the zeros $z_n$, $|\Im z_n|\simeq d_n$,
we have $\varphi_f'(x)=\sum_n|\Im z_n|/|x-z_n|^2\lesssim\sum_nd_n/((x-t_n)^2+d_n^2)$.
The terms with $t_n\notin(|x|/2,2|x|)$ contribute $\lesssim1/|x|$, and the others contribute at most
$(\min d_n)^{-1}\sum_{|x|/2<t_n<2|x|}d_n^2/((x-t_n)^2+d_n^2)\lesssim(\min d_n)^{-1}$, where
the minimum is over the indices with $t_n>|x|/2$; this tends to $0$ because
$d_n\to\infty$ (as $\psi'=o(1)$).
\vglue2pt

We consider next statement (ii). Our plan is
to construct functions $f\in {\rm PW} $ for which
$\sum_k|f(\lambda_k)|^2/\|f\|_{{\rm PW}}^2$ is arbitrarily small.
Write $\Lambda=\Lambda(\psi)$. Since $\psi=0$ on $(-\infty,0]$, we have
$\Lambda\cap(-\infty,0]=\{0,-1,-2,\dots\}$. Moreover, $\psi'$ is bounded,
so that every interval of length $1$ contains at most $K$ points of $\Lambda$,
with $K$ depending only on $\psi$.

{\it A long window.} Let $\eta,s>0$ and $\epsilon=\eta/s$. Since $x\psi'(x)\to0$,
there is $R$ such that $\psi(x)-\psi(x_0)\le\epsilon\log(x/x_0)$ when
$x\ge x_0\ge R$. As $\psi$ is continuous and unbounded, we may choose $x_0\ge R$,
as large as we wish, and an integer $n$ such that $\psi(x_0)=n+1/2$. Then
$|\psi(x)-n-1/2|\le\eta$ on $W=[x_0,x_0e^{s}]$; since $\lambda+\psi(\lambda)\in\Z$
for $\lambda\in\Lambda$, we get
$$\dist(\lambda,\Z+1/2)\le\eta\quad\hbox{for }\lambda\in\Lambda\cap W.$$
(This is the only place where $\psi'(x)=o(1/x)$ is used: $\Lambda$ stays close to
$\Z+1/2$ over arbitrarily long logarithmic ranges.)

{\it The test function.} Let $a_1=(x_0-3)/2$ and $a_2=x_0e^s+3$, so that
$a_1\ge8$ and $a_2\ge16a_1$ when $x_0$ is large and $s\ge3$. Consider the density
$m=1+\mu$, where
$$
\mu=-\frac1{2a_1}\,{\bf 1}_{[a_1,2a_1]}-\frac1{2a_2}\,{\bf 1}_{[a_2,2a_2]}
$$
consists of two negative ``ramps'' of mass $-1/2$ each; thus $1/2\le m\le1$ and
$\int_0^xm$ equals $x$ for $x\le a_1$, $x-1/2$ for $2a_1\le x\le a_2$ and $x-1$ for
$x\ge2a_2$. Apply the construction of $A$ in the proof of Corollary 2 with
$\omega'$ replaced by $m$. By Lemma 3 of \cite{OS99}, we get a real entire function $G$,
whose zeros $\gamma_k$ are the barycentres of the intervals
$[x_k,x_{k+1}]$, $\int_{x_k}^{x_{k+1}}m=1$, such that
$|G(z)|e^{-U_m(z)}\simeq\min(1,\dist(z,\Gamma))$, with constants depending only on the
bounds $1/2$ and $1$ for $m$, and not on $a_1$ and $a_2$ (this uniformity is not stated in \cite{OS99}, but it follows
from the proof given there). Let
$f(z)=e^{-\kappa(z+1/2)}G(z+1/2)$, $\kappa=\int\mu(t)\,dt/t$, which removes the linear term
$\kappa\Re z$ of $U_m=\pi|\Im z|+U_\mu$. Then:
\begin{itemize}
\item[(a)] The zeros of $f$ include the integers $\le a_1-1$, the points of
$\Z+1/2$ in $[2a_1+3,a_2-3]$, and, to the right of $2a_2+3$, again integers; in
particular $f(-j)=0$ for $j=0,1,2,\dots$.
\item[(b)] For real $x$, $|f(x)|^2\simeq{\rm env}(x)\,\dist(x,Z(f))^2$, where
$$
{\rm env}(x)=\min(1,a_1/|x|)\min(1,a_2/|x|),
$$
because $U_\mu(x)-\kappa x$ is the integral of $\log|1-x/t|$ against $\mu(t)\,dt$, and
the integral of $\log|1-x/t|$ against a mass $-1/2$ spread uniformly over $[a,2a]$ is
$-\frac12\log\max(1,|x|/a)+O(1)$. Note that
${\rm env}(x)=a_1/x$ on $[a_1,a_2]$.
\item[(c)] $U_\mu(z)-\kappa\Re z\le\log(2a_2)+1$, so that $f$ is of exponential type $\pi$, and
$f\in{\rm PW}$ by (b), since $\int{\rm env}<\infty$.
\end{itemize}

{\it Mass versus sampling sum.} On $W\subset[2a_1+3,a_2-3]$ we have
${\rm env}(x)=a_1/x$ and $Z(f)=\Z+1/2$ near $x$, and the mean value of
$\dist(x,\Z+1/2)^2$ is $1/12$; since $\int_W dx/x=s$,
\begin{equation}\label{infty}
\|f\|^2_{\rm PW}\ge\int_W|f(x)|^2dx\gtrsim a_1s.
\end{equation}
Now estimate $\sum_{\lambda\in\Lambda}|f(\lambda)|^2$. The points
$\lambda\le0$ contribute nothing, by (a). There are $\lesssim a_1$ points in
$(0,x_0)$, with $|f|^2\lesssim1$ by (b). For $\lambda>x_0e^s$ we have
$|f(\lambda)|^2\lesssim{\rm env}(\lambda)$, which is $\le a_1/\lambda$ up to $a_2$ and
$a_1a_2/\lambda^2$ afterwards, and summing over $\lesssim K$ points per unit length
gives $\lesssim a_1$. Thus
\begin{equation}\label{finite}
\sum_{\lambda\in\Lambda\setminus W}|f(\lambda)|^2\lesssim a_1.
\end{equation}
On $W$, by the window estimate, $\dist(\lambda,Z(f))\le\eta$, hence $|f(\lambda)|^2\lesssim\eta^2a_1/\lambda$, and
\begin{equation}\label{sampsum}
\sum_{\lambda\in\Lambda\cap W}|f(\lambda)|^2\lesssim\eta^2a_1\,(s+1).
\end{equation}
By \eqref{infty}, \eqref{finite} and \eqref{sampsum},
$\sum_{\lambda\in\Lambda}|f(\lambda)|^2/\|f\|^2_{\rm PW}\lesssim(1+\eta^2(s+1))/s\lesssim1/s+\eta^2$,
which is as small as we please if we first take $\eta$ small and then $s$ large. Hence
$\Lambda$ is not sampling.
\enddemo

{\it Acknowledgements}. We are indebted to Yurii Lyubarskii and
Michael Sodin for enlightening discussions.

\AuthorRefNames [OCS99]

\end{document}